# A mathematical model to assess the effects of COVID-19 on the cardiocirculatory system

Andrea Tonini[1*], Christian Vergara[2], Francesco Regazzoni[1], Luca Dedè[1], Roberto Scrofani[3], Chiara Cogliati[4] & Alfio Quarteroni[1,5]

[1]MOX, Dipartimento di Matematica, Politecnico di Milano, Milan, Italy.

[2]LABS, Dipartimento di Chimica, Materiali e Ingegneria Chimica, Politecnico di Milano, Milan, Italy.

[3]UOC Cardiochirurgia Fondazione IRCCS Ca'Granda, Ospedale Maggiore Policlinico di Milano, Milan, Italy.

[4]Ospedale L. Sacco, Milan, Italy.

[5](Professor Emeritus) Institute of Mathematics, École Polytechnique Fédérale de Lausanne, Switzerland.

*email: andrea.tonini@polimi.it

## Abstract

Impaired cardiac function has been described as a frequent complication of COVID-19-related pneumonia. To investigate possible underlying mechanisms, we represented the cardiovascular system by means of a lumped-parameter 0D mathematical model. The model was calibrated using clinical data, recorded in 58 patients hospitalized for COVID-19-related pneumonia, to make it patient-specific and to compute model outputs of clinical interest related to the cardiocirculatory system. We assessed, for each patient with a successful calibration, the statistical reliability of model outputs estimating the uncertainty intervals. Then, we performed a statistical analysis to compare healthy ranges and mean values (over patients) of reliable model outputs to determine which were significantly altered in COVID-19-related pneumonia. Our results showed significant increases in right ventricular systolic pressure, diastolic and mean pulmonary arterial pressure, and capillary wedge pressure. Instead, physical quantities related to the systemic circulation were not significantly altered. Remarkably, statistical analyses made on raw clinical data, without the support of a mathematical model, were unable to detect the effects of COVID-19-related pneumonia, thus suggesting that the use of a calibrated 0D mathematical model to describe the cardiocirculatory system is an effective tool to investigate the impairments of the cardiocirculatory system associated with COVID-19.

## Introduction

The coronavirus disease 2019 (COVID-19) caused by severe acute respiratory syndrome coronavirus 2 (SARSCoV-2) primarily affects the respiratory system, even if it does not spare other organs as it occurs for the cardiovascular system at large [1, 2, 3]. In severe COVID-19-related pneumonia, impairment of heart function seems to be mainly driven by right ventricle involvement, while consequences on the left ventricle appear to be less common [4]. Right ventricle dilation, diminished right ventricular function and elevated pulmonary arterial systolic pressure have been described and are associated with mortality in severe COVID-19 [5, 6]. Respiratory failure with shortening of oxygen supply represents the main clinical picture of the disease. Hypoxemia is associated with a huge increase in intrapulmonary shunt due to alveolar fluid filling/consolidations. In fact, the pulmonary shunt fraction (measuring the percentage of blood that does not oxygenate in the lungs) is in physiological conditions below 5% [7], whereas it reaches values up to 60% in patients with ongoing COVID-19 infection [8, 9]. Endothelial damage with diffuse micro-thrombosis has been widely described in histological studies in COVID-19 pneumonia patients and is associated with an increase in dead space in lungs and thus in non-oxygenated blood [10]. On the other hand, such an increase in intrapulmonary shunt has been postulated to depend also on an impairment of hypoxic pulmonary vasoconstriction that should restrict pulmonary flow to hypo-ventilated lung areas [11]. These mechanisms do not seem to be correlated with each other and seem to coexist to varying degrees in COVID-19 pneumonia patients [12].

In this context, physics-based mathematical models are an effective and accurate tool for making predictions through virtual scenarios and for providing clinical answers in terms of impairments of the cardiovascular function associated with COVID-19 [13, 14, 15]. In this respect, we previously studied, by means of a computational lumped-parameter (i.e., 0D) model, possible effects in terms of, e.g., cardiac output and pressures [16]. However, this previous study did not integrate clinical data into the analysis in a systematic manner.

The main novelty of this paper is to assimilate, by means of a calibration method, clinical data coming from measurements on COVID-19 patients regarding, e.g., cardiac volumes and vascular pressures, into the computational model proposed in [17] to make it patient-specific and then to use such calibrated model for making predictions on the impairments of the cardiovascular function associated with the ongoing infection. To do this, we first substituted the 3D left ventricle with a 0D component as in [16] and we improved the model of [16] by adding further compartments representing systemic and pulmonary micro-vasculatures. In this work, we focused on reproducing the blunted hypoxic pulmonary vasoconstriction, that is among the causes of the reduction in blood oxygenation, by means of the calibration of the model, thus neglecting the possible increase in pulmonary resistance associated to diffuse micro-thrombosis.

Our final goal is to study possible associations between the ongoing infection of COVID-19 and the impairments on the cardiocirculatory system by estimating physical quantities of clinical interest not available as measured clinical data and by performing a statistical analysis on these quantities.

## Results

The modified lumped-parameter model consists of a system of ODEs that has to be numerically solved to allow the computation of different model outputs of clinical interest. We calibrated the model to fit some clinical data of patients hospitalized for severe COVID-19-related pneumonia in the Internal Medicine ward of Ospedale Luigi Sacco in Milan Italy between March and April 2020.

The dataset consists of 58 patients, 29 of which were represented by the model after the calibration ($56 \pm 18$ years). All the patients required oxygen supplementation but none of them was on mechanical ventilation. Patients did not present symptoms or signs of heart failure or substantial structural cardiac disease; 10 out of 29 were older than 64 years; 6 patients had arterial hypertension, 1 had diabetes and 4 showed the association of hypertension and diabetes.

The echocardiography of each patient was performed early after the admission to the hospital. Examinations were performed at bedside using a Philips CX-50 portable device by expert operators. Measures were defined according to the latest European and American Echocardiography Society guidelines [18, 19]. We identified four groups of quantities, taken from the dataset or obtained as an output of the calibrated model:

 i) the *clinical data* used for the model calibration, obtained from clinical measurements and referring to physical quantities (PQ1), as, for example, the maximal left atrial volume (LA$_{Vmax}$) and the systolic systemic pressure (SAP$_{max}$);

ii) the *inputs* of the model (heart rate HR and body surface area BSA) and of the calibration procedure (right ventricular fractional area change $RV_{FAC}$ and tricuspid annular plane systolic excursion TAPSE), provided by other clinical measurements;
iii) the parameters of the model (e.g., resistances and compliances) determined through a calibration procedure, from now on referred to as *calibrated parameters*;
iv) the outputs of the numerical simulation of the model (e.g., flow rate and mean pressure), from now on referred to as *model outputs*. Some of them (MO1) referred to physical quantities (PQ1) that were also measured (clinical data), for example, $LA_{Vmax}$ and $SAP_{max}$. Other model outputs (MO2) referred to physical quantities (PQ2) that were not measured but quantified only by means of the computational model. Examples of the latter are the mean left atrial pressure ($LA_{Pmean}$) and indexed right ventricular end diastolic volume ($RV_{I-EDV}$). The complete list of PQ1 and PQ2 is reported in Table 1 and Table 2, respectively.

We remark that the indexed value of volumes of a patient can be computed dividing the volumes by the BSA of that patient. In what follows, an "I-" that precedes a subscript of a volume means that the volume is indexed (for example, $LV_{I-EDV}$ is the indexed left ventricular end diastolic volume).

The statistical reliability of the model outputs for which we did not have at disposal any measurement (i.e., MO2) was analysed for each patient through an estimation of the uncertainty interval, resulting from the measurement errors on clinical data. If, for a patient, the estimation of a certain model output was found to be not reliable, then this output was not used for the subsequent analysis (see *Methods – Uncertainty intervals* for further details).

We performed a statistical analysis on clinical data or model outputs MO2 (referring to physical quantities PQ1 and PQ2, respectively) to check whether their mean was significantly increased or decreased with respect to the corresponding healthy range [7, 18, 19, 20] (see Table 1 for PQ1 and Table 2 for PQ2), to highlight the impairments of the cardiocirculatory system associated with COVID-19 pneumonia. If the sample mean, calculated over all patients, of a given physical quantity fell inside the healthy range (see Table 1 for PQ1 and Table 2 for PQ2) we did not consider the physical quantity altered in association with COVID-19 infection, otherwise we performed z-tests on the sample mean. If the sample mean was less than the lower bound of the healthy range, the null hypothesis was that the mean was greater or equal than the lower bound of the healthy range, whereas the alternative hypothesis was that the mean was smaller than the lower bound of the healthy range. If we accepted the null hypothesis, then the corresponding physical quantity was considered not altered in association with the infection of COVID-19; otherwise, we considered the physical quantity altered in association with COVID-19. If, instead, the sample mean was greater than the upper bound of the healthy range, we proceeded similarly.

Notice that for group PQ1 the statistical analysis was carried out directly using the clinical data and not the MO1 values. Accordingly, the clinical data were used in a twofold way:

i) to statistically compare PQ1 clinical measures with healthy ranges independently of the application of the proposed lumped-parameter model (test I);
ii) to calibrate the lumped-parameter model for the patients at hand thus allowing to obtain MO2 that are statistically compared with healthy ranges (test II).

For the sake of clarity, we reported in Figure 1 the diagram flowchart of the followed procedure.

**Time transients of model outputs.** To perform a qualitative analysis, in Figure 2 we reported the time-dependent

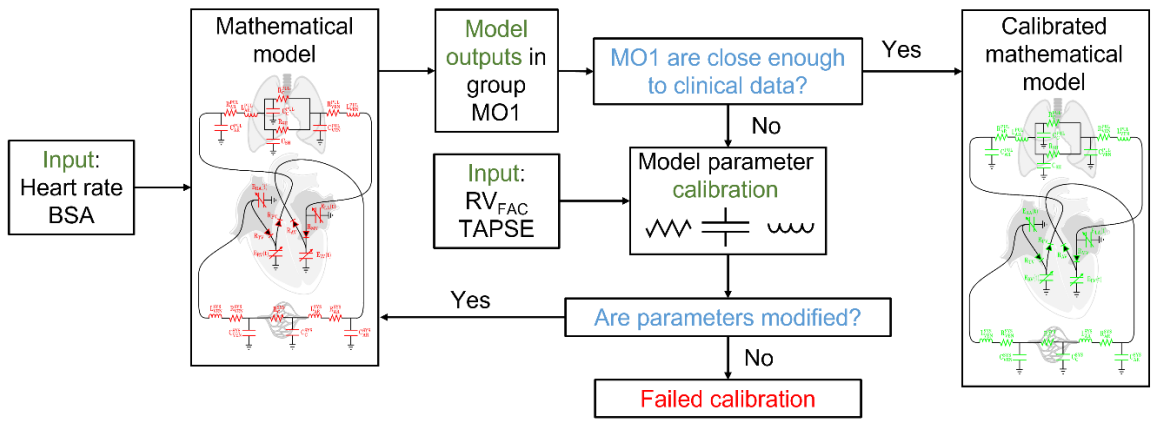
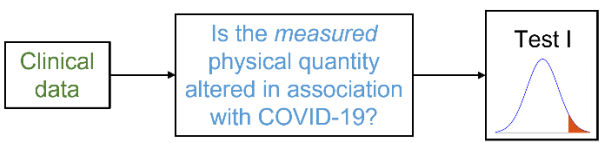
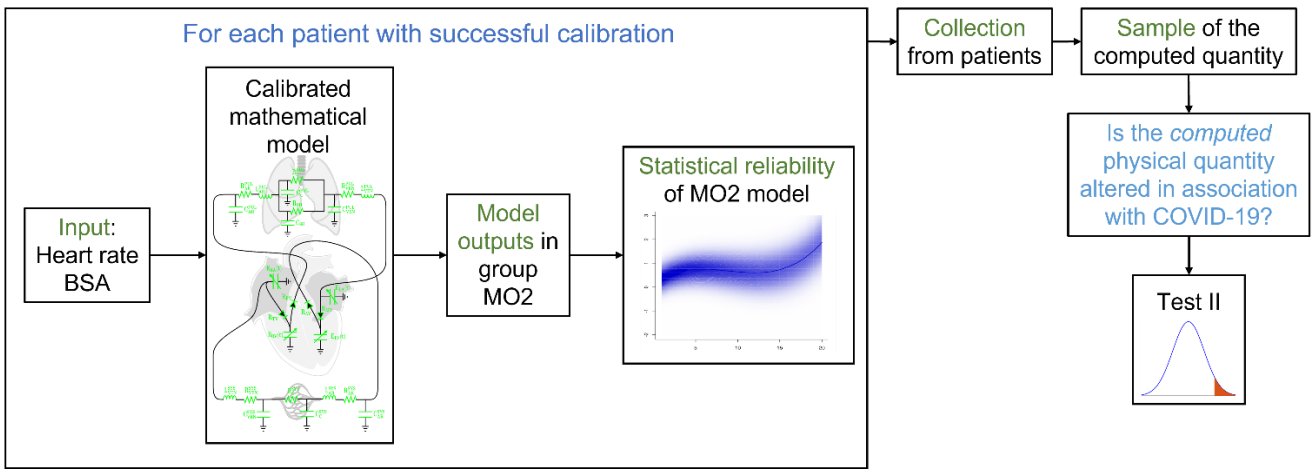

**Figure 1.** Diagram flowchart of the procedure used in this study. Top: calibration; mid: statistical analysis of measured physical quantities; bottom: statistical analysis of computed physical quantities.

Calibration: the mathematical model required as inputs HR and BSA of a specific patient. The model computed MO1 using an initial setting of parameters (that could need to be calibrated, so they are highlighted in red). If MO1 were close enough to the clinical data the model was considered calibrated (the parameters are highlighted in green); if not, the calibration method was iteratively applied to the parameters using $RV_{FAC}$ and TAPSE as inputs. If the parameters were not modified the calibration failed; if not, MO1 were recomputed by using the new setting of parameters and the previous steps were repeated.

Statistical analysis 1: we performed test I on clinical data.

Statistical analysis 2: HR and BSA were used as inputs of the calibrated model for every patient with a successful calibration, the model computed the MO2 and we checked the statistical reliability of MO2. We collected the reliable MO2 from every patient and we performed test II on the reliable MO2 of all the patients.

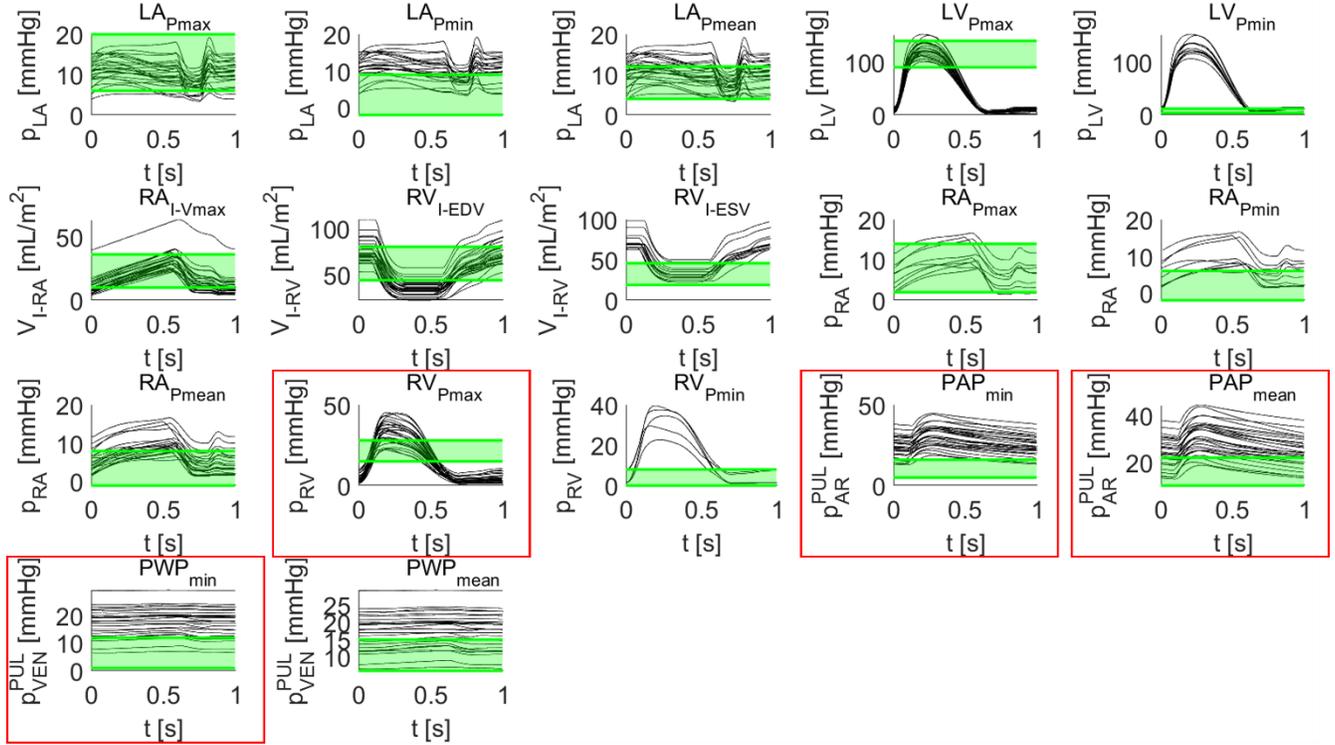

**Figure 2.** Time transients of model outputs during a cardiac cycle. In green, the reference healthy ranges of the corresponding PQ2 (the name reported above each graph) are highlighted. In red boxes the model outputs that possibly lie significantly outside of the healthy range are reported. The duration of a heartbeat was normalized to 1 s. Notice the different sample size of the plots depending on the corresponding discarded patients. The model outputs plotted are the left atrial and ventricular pressures ($p_{LA}$ and $p_{LV}$), the indexed right atrial and ventricular volumes ($V_{I-RA}$ and $V_{I-RV}$), the right atrial and ventricular pressures ($p_{RA}$ and $p_{RV}$), the pulmonary arterial and venous pressures ($p_{AR}^{PUL}$ and $p_{VEN}^{PUL}$).

model outputs (by normalizing the heartbeat duration) together with the healthy ranges (in green) related to the following physical quantities among PQ2: maximal, minimal and mean left atrial pressures ($LA_{Pmax}$, $LA_{Pmin}$ and $LA_{Pmean}$), maximal and minimal left ventricular pressures ($LV_{Pmax}$ and $LV_{Pmin}$), indexed maximal right atrial volume ($RA_{I-Vmax}$), indexed right ventricular end diastolic and systolic volumes ($RV_{I-EDV}$ and $RV_{I-ESV}$), maximal, minimal and mean right atrial pressures ($RA_{Pmax}$, $RA_{Pmin}$ and $RA_{Pmean}$), maximal and minimal right ventricular pressures ($RV_{Pmax}$ and $RV_{Pmin}$), minimal and mean pulmonary arterial pressures ($PAP_{min}$ and $PAP_{mean}$) and the minimal and mean pulmonary wedge capillary pressures ($PWP_{min}$ and $PWP_{mean}$). Notice that, for each graph, only patients such that the corresponding model output had been found to be statistically reliable (on the basis of the estimated uncertainty interval, see *Methods-Uncertainty intervals* section) were reported. We point out that, from Figure 2, some sample sizes were too small to analyse the corresponding model output (e.g., $RV_{Pmin}$). For the remaining model outputs, we moved on to the statistical analysis to study the impairments of the cardiocirculatory system associated with COVID-19, as detailed in the next paragraph.

**Statistical analysis of PQ1 (clinical data, test I) and PQ2 (test II).** We analysed the PQ1 using the clinical data (test I - Table 1), namely: the maximal left atrial volume ($LA_{Vmax}$), the left ventricular end diastolic and end systolic volumes ($LV_{EDV}$ and $LV_{ESV}$), the left ventricular ejection fraction ($LV_{EF}$), the maximal right atrioventricular pressure gradient (max $\nabla p_{rAV}$), the systolic and diastolic systemic pressures ($SAP_{max}$ and $SAP_{min}$) and the systolic pulmonary pressure ($PAP_{max}$). There was no statistical evidence that the clinical data related to PQ1 were altered in asso-ciation with COVID-19-related pneumonia because the mean of the samples lied in the corresponding healthy ranges.

| PQ1 | Healthy range | Mean ± std dev |
|---|---|---|
| $LA_{I\text{-}Vmax}$ [mL/m$^2$] ($LA_{Vmax}$ [mL]) | [16,34] [18] | 32.7 ± 13.7 (n = 56) (59.1 ± 26.2) |
| $LV_{I\text{-}EDV}$ [mL/m$^2$] ($LV_{EDV}$ [mL]) | [50,90] [20] | 56.5 ± 11.6 (n = 58) (101.5 ± 24.1) |
| $LV_{ESV}$ [mL] | [18,52] [18] | 36.8 ± 15.3 (n = 57) |
| $LV_{EF}$ [%] | [53,73] [18] | 64.5 ± 7.5 (n = 58) |
| max $\nabla p_{rAV}$ [mmHg] | - | 23.0 ± 5.9 (n = 42) |
| $SAP_{max}$ [mmHg] | [-,140] [18] | 120.6 ± 14.7 (n = 58) |
| $SAP_{min}$ [mmHg] | [-,80] [18] | 71.0 ± 11.4 (n = 58) |
| $PAP_{max}$ [mmHg] | [15,28] [20] | 27.9 ± 5.1 (n = 40) |

**Table 1.** Statistics of clinical data. The mean and standard deviation of the samples are provided together with the sizes in brackets. Notice that the sample sizes of the clinical data are different due to heterogeneous samples. The hypothesis tests were not performed because the mean of the samples lied in the respective healthy range. Test I.

Instead, we analysed the PQ2 using the model outputs MO2 (test II - Table 2), obtaining the following outcomes:

I. For $RV_{Pmax}$, $PAP_{min}$, $PAP_{mean}$, $PWP_{min}$ and $PWP_{mean}$ we rejected the null hypothesis and thus these physical quantities resulted significantly increased with respect to the healthy ranges;
II. For left ventricular stroke volume ($LV_{SV}$), cardiac index (CI) and thus the cardiac output (CO), $LA_{Pmax}$, $LA_{Pmean}$, $LV_{Pmax}$, $RA_{I\text{-}Vmax}$, $RV_{I\text{-}EDV}$, right ventricular ejection fraction ($RV_{EF}$), systemic and pulmonary vascular resistances SVR and PVR, we did not reject the null hypothesis, thus there was no statistical evidence that these physical quantities were altered in association with COVID-19-related pneumonia;
III. The sample sizes of $LA_{Pmin}$, $LV_{Pmin}$, $RV_{I\text{-}ESV}$, $RA_{Pmax}$, $RA_{Pmin}$, $RA_{Pmean}$, $RV_{Pmin}$ and the Shunt Fraction were too small to perform the hypothesis tests.

## Discussion

This study addressed the association between COVID-19-related pneumonia and the impairments of the cardiovascular system. This has been faced by analysing clinical measures and model outputs computed through a calibrated lumped-parameter cardiocirculatory mathematical model. To the best of our knowledge, the current study is the first that used clinical measures and calibrated models to infer the cardiovascular physical quantities significantly altered in association with COVID-19-related pneumonia.

We start by discussing the available clinical data measured at Ospedale Luigi Sacco in Milan and related to cardiovascular physical quantities for COVID-19 pneumonia patients. We found that none of the measured physical quantities (i.e., PQ1) was altered in association with COVID-19-related pneumonia (Table 1). See also [21] for another analysis of the same dataset of clinical measures.

Regarding the analysis of MO2, we noticed from Figure 2 that some of the related physical quantities (among PQ2) lied within healthy ranges (e.g., $LV_{Pmax}$), whereas other physical quantities lied outside them (e.g., $RV_{Pmax}$). For the remaining physical quantities, we could not infer from Figure 2 if they were altered or not in association with COVID-19-related pneumonia (e.g., $PWP_{mean}$ or $RA_{Pmean}$). Therefore, to significantly assess the alterations associated with COVID-19, we resorted to hypothesis tests.

|  | PQ2 | Healthy range | Mean ± std dev (HP1) | test II (p-value) | COVID-19 literature values |
|---|---|---|---|---|---|
| **I) Rejected null hypothesis** | $RV_{Pmax}$ [mmHg] | [15,28] [20] | 33.7 ± 6.8 (n = 29) | **2.62E-06** | [30,46] [22] |
|  | $PAP_{min}$ [mmHg] | [5,16] [20] | 23.6 ± 6.2 (n = 29) | **3.06E-11** | [15,26] [11] |
|  | $PAP_{mean}$ [mmHg] | [10,22] [20] | 27.1 ± 6.5 (n = 29) | **9.97E-06** | [25,33] [11] |
|  | $PWP_{min}$ [mmHg] | [1,12] [20] | 17.1 ± 5.2 (n = 28) | **1.04E-07** | - |
|  | $PWP_{mean}$ [mmHg] | [6,15] [20] | 17.5 ± 5.1 (n = 28) | **5.35E-03** | [11,18] [11] |
| **II) Not rejected null hypothesis** | $LV_{SV}$ [mL] | [30,80] [18] | 74.0 ± 10.7 (n = 29) | - | [68,105] [11] |
|  | CI [L/min/m$^2$] (CO [L/min]) | [2.8,4.2] [20] | 3.2 ± 0.5 (n = 29) (5.9 ± 1.0) | - | [2.7,4.5] [11] / [1.98,3.32] [23] ([4.4,6.3] [24]) |
|  | $LA_{Pmax}$ [mmHg] | [6,20] [20] | 12.8 ± 3.2 (n = 25) | - | - |
|  | $LA_{Pmean}$ [mmHg] | [4,12] [20] | 10.2 ± 2.8 (n = 27) | - | - |
|  | $LV_{Pmax}$ [mmHg] | [90,140] [20] | 124.4 ± 13.3 (n = 29) | - | - |
|  | $RA_{I-Vmax}$ [mL/m$^2$] | [10,36] [18] | 31.8 ± 8.0 (n = 28) | - | [15,29] [11] / [14,25] [23] |
|  | $RV_{I-EDV}$ [mL/m$^2$] | [44,80] [19] | 75.4 ± 12.4 (n = 29) | - | - |
|  | $RV_{EF}$ [%] | [44,71] [19] | 53.6 ± 5.3 (n = 29) | - | - |
|  | SVR [mmHg min/L] | [11.3,17.5] [20] | 15.9 ± 3.3 (n = 29) | - | [8.1,13.0] [11] |
|  | PVR [mmHg min/L] | [1.9,3.1] [20] | 3.0 ± 1.4 (n = 28) | - | [3.1,4.7] [11] |
| **III) Sample size too small** | $LA_{Pmin}$ [mmHg] | [-2,9] [20] | 7.4 ± 2.2 (n = 22) | - | - |
|  | $LV_{Pmin}$ [mmHg] | [4,12] [20] | 6.2 ± 1.4 (n = 12) | - | - |
|  | $RV_{I-ESV}$ [mL/m$^2$] | [19,46] [19] | 33.1 ± 8.4 (n = 14) | - | - |
|  | $RA_{Pmax}$ [mmHg] | [2,14] [20] | 11.7 ± 3.3 (n = 9) | - | - |
|  | $RA_{Pmin}$ [mmHg] | [-2,6] [20] | 4.4 ± 2.9 (n = 10) | - | - |
|  | $RA_{Pmean}$ [mmHg] | [-1,8] [20] | 7.0 ± 2.6 (n = 23) | - | - |
|  | $RV_{Pmin}$ [mmHg] | [0,8] [20] | 3.0 ± 2.6 (n = 5) | - | - |
|  | Shunt Fraction [%] | [0,5] [7] | 3.7 ± 0.8 (n = 9) | - | - |

**Table 2.** Statistics of MO2. The mean and the standard deviation of samples are provided together with the sizes in brackets. If there is statistical evidence of impairments of the cardiocirculatory system associated with COVID-19 the p-value is highlighted in red. If the mean of a sample lied in the healthy range, the hypothesis test was not performed. The sample sizes too small to perform the hypothesis tests are highlighted in orange. For some of the physical quantities with a big sample size, we report the COVID-19 values taken from literature. Test II.

We found that the pulmonary resistances (PVR), did not significantly increase in association with COVID-19-related pneumonia (Table 2). Nonetheless, we highlighted a slightly large value of $PAP_{max}$ (Table 1) that was accompanied by a significant increase not only in $PAP_{min}$, $PAP_{mean}$ and $RV_{Pmax}$, but also in $PWP_{min}$ and $PWP_{mean}$ (Table 2). These results seem to be in line with previous evidence reported in COVID-19-related pneumonia patients studied with cardiac catheterization [11]. In this study, patients did not show an increase in PVR but the mild increase in pulmonary arterial pressure was associated with an increase in wedge pressure. The authors hypothesized that a hyperdynamic state not accompanied by an increased in hypoxic-driven vasoconstriction could determine (especially in their population of old and often hypertensive patients) an increase in wedge pressure related to an increase of LV filling pressure. In our population a substantial percentage of patients were old (34% were older than 64 years), with arterial hypertension and/or diabetes, conditions that could be in line with this interpretation, taking into account that the mean value of cardiac output computed by the model was rather large (5.9 ± 1.0 L/min, Table 2).

There was no statistical evidence that the maximal and mean left atrial pressures increased (Table 2). This could be due to limitations of the lumped-parameter model in representing the atria. Unfortunately, the sample size of $LV_{Pmin}$ was too small to infer any interpretation.

In what follows, we refer to clinical literature of patients affected by COVID-19 for a comparison with the outcomes of our mathematical model (the model outputs MO2) [11, 22, 23, 24] (see Table 2). If the mean of our samples lied in the intervals identified in clinical literature, we considered them in accordance one another. We noticed from Table 2 that the sample mean of some of the physical quantities ($RV_{Pmax}$, $PAP_{min}$, $PAP_{mean}$, $PWP_{mean}$, $LV_{SV}$ and CO) agreed with the COVID-19 literature, whereas the means of $RA_{I-Vmax}$ and SVR were slightly larger and PVR slightly lower than the values of literature, although still lying inside the healthy range.

We emphasise that the statistical analysis of raw clinical data did not allow us to infer alterations in the cardiovascular system in association with COVID-19 infection (Table 1). Instead, thanks to the computational model we proposed, suitably calibrated by using the clinical data, we were able to identify some specific physical quantities related to pulmonary circulation (i.e., $RV_{Pmax}$, $PAP_{min}$, $PAP_{mean}$, $PWP_{min}$ and $PWP_{mean}$) which were significantly altered in association with COVID-19, in the sense that there was a statistically relevant discrepancy with respect to the healthy ranges. This showed the importance of combining clinical data and computational models as an effective strategy to give meaningful insights about the impairments of the cardiocirculatory system associated with COVID-19 on cardiovascular physical quantities, which was not possible with raw clinical data and non-calibrated computational tools.

We now discuss the limitations of this study. First, notice that we did not have at disposal a control group to perform hypothesis tests in tests I-II, so we took a conservative approach comparing the mean of our samples with the lower and upper bounds of healthy ranges found in literature to infer the impairments of the cardiovascular system in association with the infection of COVID-19. The sample means of physical quantities significantly outside the corresponding healthy range highlight a clear impairment of a compartment of the cardiocirculatory system in association with COVID-19 infection. Nevertheless, we do not exclude that small changes in some physical quantity could indicate an impairment in the cardiocirculatory system as well.

Second, although being able to capture the considerable haemodynamic features, the lumped-parameter model is rather simple in comparison to other models for the study of the cardiac function (see e.g. [25, 26, 27, 28]). Improvements of the computational model will allow also to use other clinical measurements not used in this work (such as those based on partial pressures of oxygen and carbon dioxide).

Third, to quantify the uncertainty in the estimation of the model outputs, we adopted a rather simple approach in terms of independence in the selection of the parameter configurations used to computationally generate the model outcomes. More sophisticated strategies that account for a selective choice of the new parameter configuration starting from the previous ones (e.g., Markov chain Monte Carlo methods [29]), which nevertheless entail a larger computational cost, could be considered in further developments of this work.

Possible improvements of the present work are also related to the clinical measurement acquisition. Other clinical measurements, when available, could be added to the framework of the present work to improve the outcomes. It may be of particular interest having a measure of the shunt fraction, that gives information on the pulmonary capillaries, to avoid the a priori assumption between micro-thrombosis and blunted hypoxic pulmonary vasoconstriction. As a limitation, in this work we neglected the contribution of micro-thrombosis [10] and we focused only on the study of blunted hypoxic pulmonary vasoconstriction in the increase of non-oxygenated blood [11].

# Methods

**Mathematical model.** The cardiovascular system was studied by means of a lumped-parameter (0D) mathematical model that splits the system into compartments (e.g. right atrium, systemic arteries/veins) and, for each of them, the time evolution of model outputs (pressures, flow rates and cardiac volumes) is modelled by a system of ODEs [30, 31]. The lumped-parameter model is described through an electrical circuit analogy: the current represents the blood flow through vessels and valves; the electric potential the blood pressure; the electric resistance plays the role of the resistance to blood flow; the capacitance represents the vessel compliance; the inductance the blood inertia; the increase of elastance the cardiac contractility.

There are different possible choices and number of compartments, depending on the purpose of the study, for the construction of a lumped-parameter model (e. g. [16, 17, 32, 33]). We considered the computational model introduced in [17], wherein the four heart chambers, the systemic and pulmonary circulations, with their arterial and venous compartments were included, and we substituted the 3D left ventricle with a 0D component and we added two new compartments accounting for systemic and pulmonary capillaries. The pulmonary capillary circulation was also split in two compartments accounting for oxygenated and non-oxygenated capillaries (Figure 3).

The system of ODEs associated with the lumped-parameter model is formed by the equations representing continuity of flow rates at nodes and of pressures in the compartments, and its numerical solution allows to compute several model outputs as functions of time: the left and right atrial and ventricular volumes ($V_{LA}$, $V_{LV}$, $V_{RA}$ and $V_{RV}$), the systemic and pulmonary arterial, capillary and venous pressures ($p_{AR}^{SYS}$, $p_{C}^{SYS}$, $p_{VEN}^{SYS}$, $p_{AR}^{PUL}$, $p_{C}^{PUL}$ and $p_{VEN}^{PUL}$), the systemic and pulmonary arterial and venous blood fluxes ($Q_{AR}^{SYS}$, $Q_{VEN}^{SYS}$, $Q_{AR}^{PUL}$ and $Q_{VEN}^{PUL}$).

Starting from these functions, it is possible to compute the pressures of the four cardiac chambers ($p_{LA}$, $p_{LV}$, $p_{RA}$ and $p_{RV}$), the blood fluxes through the valves ($Q_{MV}$, $Q_{AV}$, $Q_{TV}$ and $Q_{PV}$), through the systemic capillaries ($Q_{C}^{SYS}$) and through oxygenated and non-oxygenated pulmonary capillaries ($Q_{C}^{PUL}$ and $Q_{SH}$), and all the model outputs referring to PQ1 and PQ2.

We considered reference values of the parameters (such as resistances and compliances) such that all the model outputs were in the reference healthy ranges of the corresponding physical quantities taken from the literature [7, 18, 19, 20] for an ideal individual with HR equal to 80 bpm (beats per minute) and BSA equal to 1.79 m². We did not consider model outputs computed starting from the flow rates, because they are not uniquely defined depending on the tract of the compartment where they are measured, from $p_{C}^{SYS}$, due to the heterogeneity of the pressures of systemic capillaries among tissues, and from $p_{VEN}^{SYS}$, even if we recovered the value of central venous pressure, that coincides with the right atrial pressure [20].

The lumped-parameter model was numerically discretized by means of Dormand-Prince method [34] (adaptive stepsize Runge-Kutta) which was implemented in Python using the Jax library [35].

**Calibration.** The lumped-parameter model was characterized by parameters representing the functional properties of the compartments (e.g., resistances). To properly select such values for a specific compartment and patient, a calibration procedure was needed [29, 36].

The calibration of the model relies on the method we presented in [29], that is aimed to reduce the sum of squared relative errors between the model outputs MO1 and clinical data, modifying the parameters of the model in suitable bounded intervals $I_i$, for $i = 1, ..., N_p$, where $N_p$ is the number of parameters independent of the patient, built starting from the reference values of parameters mentioned before.

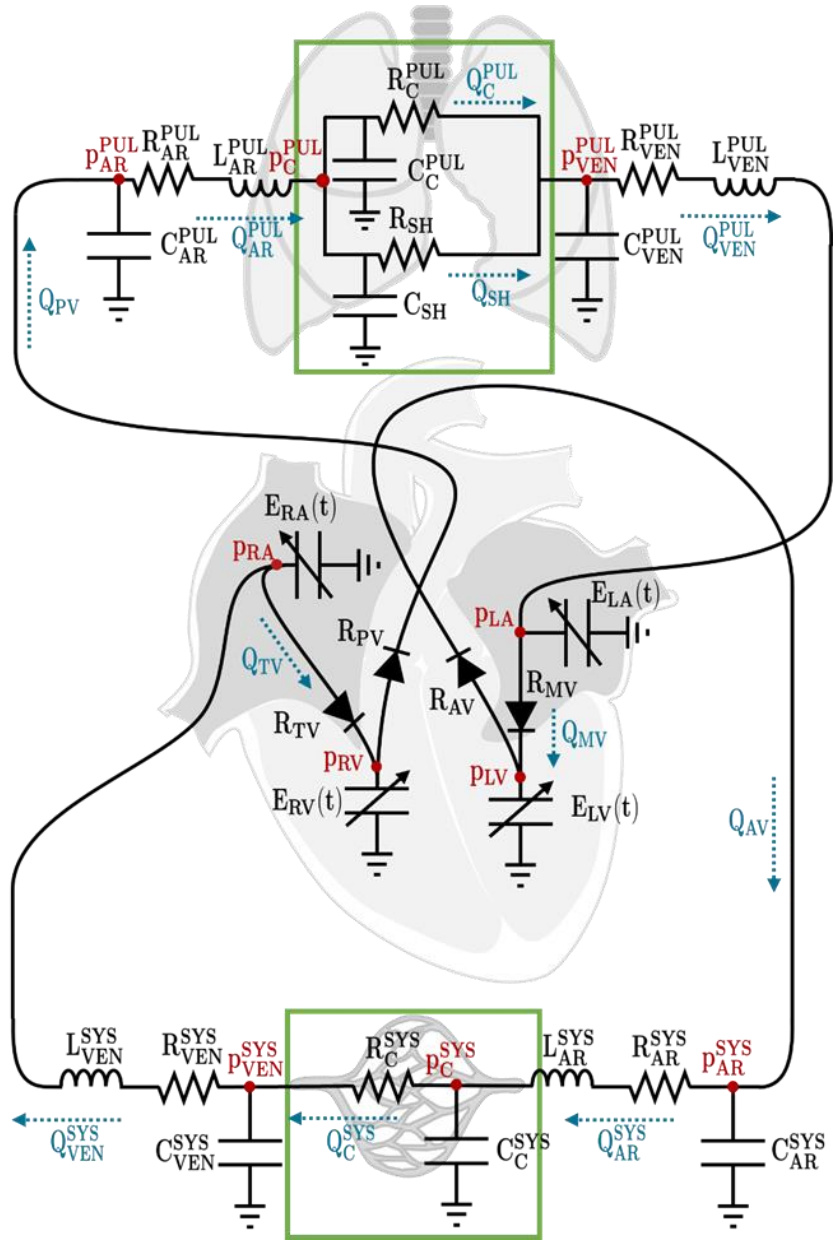

**Figure 3.** Lumped-parameter cardiocirculatory model. The unknown pressures and flow rates are in red and blue, respectively, whereas the model parameters are in black. Notice in the green boxes the new compartments with respect to [16] featuring this work.

The calibration was based on clinical measurements of COVID-19 patients that were provided by Ospedale Luigi Sacco in Milan and referred to HR and BSA, which were used as inputs for the lumped-parameter model, RV$_{FAC}$ and TAPSE, which determined the bounded interval $I_{\bar{i}}$ used during the calibration, with $\bar{i}$ the index referring to the right ventricular active elastance, and the clinical data, given by a subset of the pressures and volumes involved in the cardiac circulation.

To provide further mathematical details, we indicate with **p** a configuration of parameters of the cardiocirculatory model. The calibration method aimed to find the configuration of parameters $\bar{\mathbf{p}}^j$ which minimized the loss function for the specific patient $j$, that reads:

$$L^j(\mathbf{p}) = \sum_{l=1}^{N^j} \left( \frac{q_{m_j(l)}^j(\mathbf{p}) - d_l^j}{d_l^j} \right)^2, \quad (1)$$

where $N^j$ is the number of available clinical data for patient $j$, $d_l^j$ is the value of the l-th clinical data of patient $j$ and $q_{m_j(l)}^j$ is the value of the model output related to the l-th clinical data of patient $j$. The index $m$ of $q_m^j$ lies in $\{1, \ldots, N_q\}$ where $N_q$ is the number of both MO1 and MO2. We considered the model calibrated for a specific patient if the loss function was below $10^{-3}$. Notice that, for some patients, the calibration procedure could fail, if, for example, it reaches the minimum of the loss function that is above the required threshold.

During the calibration procedure the parameters $\mathbf{p}$ could vary in suitable intervals. Notice that, to reproduce the hypoxic pulmonary vasoconstriction condition, the resistance of non-oxygenated pulmonary capillaries ($R_{SH}$) could decrease in such a way the shunt fraction could reach values up to 70% in the worst-case scenario.

Moreover, to improve the robustness of the calibration procedure, we repeated, for every patient, the calibration three times, with different initial configurations of parameters, and we considered the calibrated setting of parameters that returned the lowest loss function.

The loss function (1) was minimized by the Quasi-Newton method L-BFGS-B [37] implemented in Scipy by computing its gradient by means of automatic differentiation (reverse mode gradient) included in the library Jax [35].

**Uncertainty intervals.** For every patient $j$ calibrated with a loss function below $10^{-3}$, a configuration of parameters $\bar{\mathbf{p}}^j$ was at disposal. The loss function was computed using the clinical data provided by Ospedale Luigi Sacco, which were related to *measurement errors*, that also affected the uncertainty of the model outputs $\mathbf{q}^j$. We needed to determine, for every patient, if the related model outputs were reliable or not, so we proceeded along two steps:

1. Build a sample of candidate model outputs $\mathbf{q}^{j,k}$ for $k = 1, \ldots, n$ ($n$ was 100);
2. Determine, by employing a simple statistical analysis, whether the mean of the model outputs was reliable.

Regarding step 1, for every provided clinical data $d_l^j$ of patient $j$, we built an interval $M_l^j$ centred in the value of the clinical data with width equal to two times the measurement error. Then, we built the samples $\mathbf{q}^{j,k}$ by following the subsequent procedure:

a) Choose a relative width $w$ ($w$ was 12.5%);
b) Build an interval centred at $\bar{p}_i^j$ and with width $2w\bar{p}_i^j$ for every $i = 1, \ldots, N_p$. If this interval is not included in the parameter interval $I_i$ used for the calibration, then cut off its overflowing extremities.
c) Perturb every parameter of the calibrated patient sampling from a uniform distribution in the corresponding interval built at point b) thus obtaining $p_i^j$;
d) Run a simulation of the cardiocirculatory model with parameters $\mathbf{p}^j$;
e) Check if the model output $q_{m_j(l)}^j$ generated at point d) lie in the intervals $M_l^j$. If they do, save the new configuration of acceptable model outputs $\mathbf{q}^j$, otherwise reject it;
f) Repeat from point c) until $n$ iterations are performed;

g) Check if the acceptance ratio (ratio between the number of saved configurations and the number of iterations) is within $[0.1, 0.15]$. If it does, repeat from point c) to e) until $n$ configurations are accepted because at this step the sample size of candidate model outputs is small (with $n = 100$, the size is between 10 and 15), otherwise increase or decrease $w$ to retrieve the condition on the acceptance ratio, discard the previous configurations and repeat from point b).

Once the above procedure was concluded, we proceeded with step 2 by using the $n$ samples of acceptable model outputs $\mathbf{q}^{j,k}$ for $k = 1, \ldots, n$ generated at the previous step, for every specific patient $j$. If the standard deviation of the sample of a model output of patient $j$ was lower than 5% of its mean, we considered the mean reliable and we used it for the hypothesis tests. In this way, for every model output we built a sample of accepted values (depending on the patient), where sample size depended on the considered model output.

Prediction intervals could have been used for this analysis, but, if the sample was not normally distributed, a link function would be needed to retrieve normality [38]. We checked, for every patient $j$ and for every model output, if the sample of that model output was normally distributed by means of a chi-squared test. It turned out that the sample is not normally distributed for all patients. Thus, since we wanted to use the same statistical approach for every patient, we resorted to this heuristic approach based on standard deviation instead of prediction intervals.

**Statistical analysis.** If the sample mean of a clinical data or MO2 lied outside the corresponding healthy range we performed z-tests to check whether the mean of the model output lied significantly (p-value below 0.01) outside the healthy ranges of the corresponding physical quantity [7, 18, 19, 20] to investigate the impairments of the cardiovascular system in association with COVID-19 infection.

For each clinical datum, we computed the mean and the standard deviation of its sample without resorting to the mathematical model. The sample sizes were large enough to use one-tailed z-tests (assuming the variance equal to the unbiased sample variance) comparing their means to the nearest bound of the healthy range (test I). The tests were left or right tailed if the sample mean lied at the left or at the right of the healthy range of the related physical quantity, respectively.

For every MO2 we computed the mean and the standard deviation of its sample. We performed a chi-squared test and not every sample was normally distributed, so we opted for one-tailed z-tests (assuming the variance equal to the unbiased sample variance) only if the sample had more than 25 elements comparing their means to the nearest bound of the healthy range (test II). The tests were left or right tailed if the sample mean lied at the left or at the right of the healthy range of the related physical quantity, respectively.

# Data availability

The datasets generated and/or analysed are available from the corresponding author upon reasonable request.

## Acknowledgments

This work has been supported by the Italian research project FISR (Fondo Integrativo Speciale per la Ricerca) 2020 " Mathematical modelling of Covid-19 effects on the cardiac function, Mathematical modelling and analysis of clinical data related to the Covid-19 pandemic in Italy, 2021_ASSEGNI_DMAT_10 ". Funding agency: MIUR (Italian Ministry of Education, Universities and Research).

## Author contributions

Acquisition of the clinical data: R.S., C.C.

Methodology: A.T., F.R.

Conceptualization: R.S., C.C., A.T., F.R., C.V., L.D., A.Q.

Formal analysis and investigation: A.T.

Interpretation of the results: R.S., C.C., A.T.

Writing - Original draft preparation: A.T., C.V.

Writing - Review and editing: R.S., C.C., A.T., F.R., C.V., L.D., A.Q.

Figures: A.T., F.R., C.V.

Supervision: C.V., A.Q.

## Competing interests

The authors declare no competing interests.

## Additional information

**Correspondence** and requests for materials should be addressed to andrea.tonini@polimi.it